 \documentclass[preprint,3p,times]{elsarticle}

\usepackage{amssymb}

\usepackage{graphics} 
\usepackage{graphicx}
\usepackage{amsmath} 
\usepackage{amsfonts}
\usepackage{mathrsfs}

\newcommand{\hbd}{\mathbf{H}} 
\newcommand{\B}{\mathcal{B}} 
\newcommand{\bk}{\mathfrak{B}} \newcommand{\ik}{\mathfrak{I}}
\newcommand{\ok}{\mathfrak{O}}
\newcommand{\rs}{\mathcal{R}} 
\newcommand{\xs}{\mathcal{X}}
\newcommand{\set}[1]{\left\{#1\right\}}

\newcommand{\abs}[1]{\left\vert #1 \right\vert}
\newcommand{\ra}{\rightarrow}
\newcommand{\sra}{\rightrightarrows}
\newcommand{\R}{\mathbb{R}}

\newcommand{\N}{\mathbb{N}}
\newcommand{\soln}{\mathfrak{S}}
\newcommand{\inte}{\operatorname{int}}
\newcommand{\dom}{\operatorname{dom}}
\newcommand{\rge}{\operatorname{rge}}
\newcommand{\eps}{\varepsilon}
\newcommand{\dt}{\delta}
\newcommand{\sub}{\subseteq}
\newcommand{\wh}{\widehat}
\newcommand{\xk}{\mathfrak{x}}
\newcommand{\zero}{\mathbf{0}}
\newcommand{\zk}{\mathfrak{z}}
\newcommand{\uf}{\mathfrak{u}}
\newcommand{\uu}{\mathcal{U}}
\newcommand{\ep}{\vartheta}

\newtheorem{thm}{\bf Theorem}[section]
\newtheorem{cor}[thm]{\bf Corollary}
\newtheorem{lem}[thm]{\bf Lemma}
\newtheorem{rem}[thm]{\bf Remark}
\newtheorem{ass}[thm]{\bf Assumption}
\newtheorem{prop}[thm]{\bf Proposition}
\newtheorem{prob}{\bf Problem}
\newtheorem{deff}[thm]{\bf Definition}

\newcommand{\ymmark}[1]{{\color{black} #1}}

\DeclareMathAlphabet{\mathcal}{OMS}{cmsy}{m}{n}

\newcommand{\BlackBox}{\rule{1.5ex}{1.5ex}}    
\newcommand{\pfbox}{\hfill\BlackBox\\[2mm]}  
\newenvironment{pf}{\par\noindent{\bf Proof\ }}{\hfill\BlackBox\\[2mm]}

\usepackage{xcolor}

\usepackage{hyperref}
 \hypersetup{
    colorlinks=true,
    linkcolor=blue,
    filecolor=magenta,      
    urlcolor=cyan,
    pdftitle={Overleaf Example},
    pdfpagemode=FullScreen,
    }
    \usepackage[misc]{ifsym}

\journal{Nonlinear Analysis: Hybrid Systems}

\begin{document}

\begin{frontmatter}



\title{Lyapunov-Barrier  Characterization of Robust  Reach-Avoid-Stay Specifications for Hybrid Systems\tnoteref{label1}}
\tnotetext[label1]{This paper was not presented at any IFAC meeting. }


\author{Yiming Meng\corref{cor1}}
\ead{yiming.meng@uwaterloo.ca} 
\author{Jun Liu\corref{cor1}}
\ead{j.liu@uwaterloo.ca}
\affiliation{organization={Department of Applied Mathematics, University of Waterloo},
    city={ Ontario},
    country={Canada}}

\cortext[cor1]{Corresponding  authors.}

\begin{abstract}
Stability, reachability, and safety are crucial properties of dynamical systems. While verification and control synthesis of reach-avoid-stay objectives can be effectively handled by abstraction-based formal methods, such approaches can be computationally expensive due to the use of state-space discretization. In contrast,  Lyapunov methods qualitatively characterize stability and safety properties without any state-space discretization. 
Recent work on converse
Lyapunov-barrier theorems also demonstrates an approximate completeness for verifying reach-avoid-stay specifications of systems modelled by nonlinear  differential equations. In this paper, \ymmark{based on the topology of hybrid arcs}, we extend the Lyapunov-barrier characterization to more general hybrid systems described by differential and difference inclusions. We show that Lyapunov-barrier functions are not only sufficient to guarantee reach-avoid-stay specifications for well-posed hybrid systems, but also necessary for arbitrarily slightly perturbed systems under mild conditions. \ymmark{Numerical examples are provided to illustrate the main results.}
\end{abstract}

\begin{keyword}
Hybrid systems; Reach-avoid-stay specifications;  Stability with safety; Robustness; Converse Lyapunov-barrier functions theorems.  



\end{keyword}

\end{frontmatter}


\section{Introduction}
The reach-avoid-stay property is one of the building blocks
for specifying more complex temporal logic objectives. States (or controlled states) of dynamical systems with  reach-avoid-stay properties are always within the
safe region and, moreover, reach a target set within a finite
time and stay inside it afterwards. 

For hybrid systems modelled by differential and difference equations with possibly small perturbations, verification and control synthesis of various linear temporal logic properties are achievable via abstraction-based formal
methods \cite{kloetzer2008fully, tabuada2006linear, nilsson2017augmented} and model checking algorithms \cite{baier2008principles}, where state space discretization and abstraction analysis are needed to guarantee the soundness and (approximated) completeness \cite{girard2009approximately, pola2008approximately,  liu2017robust, liu2021closing} of the generated symbolic models. Regardless of the potentially heavy computational efforts, formal methods are guaranteed to return  the set of initial states from which the specification are satisfied. The current literature that pertains to formal methods focuses on developing algorithms to reduce the computational complexity \cite{gao2014delta, girard2015safety, li2020robustly, li2022robustly}. However, case studies in \cite{meng2021control} have shown that for systems that also continuously depend on tunable parameters and undergo bifurcations, challenges arose regarding the sampling time for constructing abstractions under the non-uniform rates of state evolution.

Compared to formal methods, Lyapunov-like functions are able
to characterize stability and safety related attributes without state space discretization \cite{agrawal2017discrete, hsu2015control, nguyen2020dynamic}. The flexible quadratic programming framework \cite{ames2016control} with Lyapunov and barrier constraints also shows effectiveness
in the control synthesis of  stabilization and safety objectives. The work in \cite{romdlony2016stabilization} developed  sufficient conditions for stabilization with safety
guarantees by uniting Lyapunov and barrier functions (see \cite[Proposition 1]{romdlony2016stabilization} for details). However, the
conditions  imply that the unsafe set has 
to be unbounded \cite{braun2017existence}. The authors in \cite{braun2020comment} further made an interesting observation that the derived control Lyapunov-barrier conditions proposed in \cite[Section 4, 5]{romdlony2016stabilization} are topologically invalid. They
instead proposed sufficient conditions for safe stabilization
using non-smooth control Lyapunov functions \cite{braun2019complete}. 

The recent work 
\cite{meng2022smooth} improved the Lyapunov-barrier approach and 
showed that it is possible to construct united Lyapunov and barrier functions, or even a single Lyapunov function,  that
are defined on the entire set of initial conditions from which
stabilization with safety guarantees is satisfied. The connection between stabilization with safety guarantees and reach-avoid-stay was also established via robustness. The work \cite{meng2021control}  proposed control Lyapunov-barrier functions and
provided sufficient conditions for control
synthesis with reach-avoid-stay specifications. Such a method has been numerically shown to outperform abstraction-based algorithms in a system that undergoes a Hopf-bifurcation. 

The latest research \cite{meng2022acc} has shown that the existence of the stochastic analogue of Lyapunov-barrier functions is also sufficient for characterizing probabilistic
reach-avoid-stay specifications for stochastic systems with point-mass uncertainties.  However, due to the effects of diffusion and the concept of weak solutions, probabilistic reach-avoid-stay specifications, other than strong reach-avoid-stay with probability one, may fail to be related to probabilistic stability with safety guarantees w.r.t. some subset of the target set, and hence the (approximated) necessity of Lyapunov-barrier characterization for probabilistic reach-avoid-stay does not hold for stochastic dynamical systems.

In this paper, we show that the Lyapunov-barrier approach can be  extended to  verification of reach-avoid-stay specifications for hybrid systems with differential and difference inclusions. Unlike systems with diffusion,  reach-avoid-stay properties of solutions to hybrid systems can be converted to stabilization with safety guarantees under mild conditions. The approximate equivalence of these two types of specifications can be established via robustness. We show that smooth  Lyapunov-barrier functions can be constructed given the compactness of target set and the set of initial states. \ymmark{Note that even though the idea follows the previous work in 
\cite{meng2022smooth} and \cite{meng2022acc}, the underlying
topology of the solutions to hybrid systems are different. Hence, the previous results are not directly applicable. We aim to leverage the rich results on stability theory of hybrid systems \cite{goebel2012hybrid} to unify Lyapunov and barrier conditions in the context of converse Lyapunov-barrier theorems.}

We mention that the latest converse barrier theorems for differential inclusions  \cite{ghanbarpour2021barrier, maghenem2022converse,maghenem2022strong}, in comparison to the previous works \cite{wisniewski2015converse, ratschan2018converse, liu2021converse},  provide a possibility to construct less smooth Lyapunov-barrier functions under less restricted topological requirement, e.g., unbounded reachable set. In this paper,  we only consider the case where the uniformly asymptotically stable set is compact, which seems already satisfactory in practice. For this reason, we stick with the compactness assumption on the target set and the set of initial states, and discuss the possible relaxation of safe conditions based on \cite{ghanbarpour2021barrier, maghenem2022converse,maghenem2022strong} in Section \ref{sec: converse}.

The rest of this paper is organized as follows. In Section
\ref{sec: pre}, we present the preliminaries for the systems, concepts of
solutions, as well as other important definitions. In Section \ref{sec: connect}, we show the connections between robust 
reach-avoid-stay and stability with safety guarantees. In Section \ref{sec: converse},
we provide
the converse Lyapunov-barrier functions theorem for both stability with safety guarantees and reach-avoid-stay specifications. We illustrate the results with examples in Section \ref{sec: example} and
conclude the paper in Section \ref{sec: conclusion}.

\textbf{Notation:} We denote by $\R^n$ the Euclidean space of dimension $n>1$, by $\R$ the set of real numbers, by $\R_{>0}$ the set of positive real  numbers, and by $\R_{\ge 0}$  the set of nonnegative real  numbers. For $x\in\R^n$ and $r\ge 0$, we denote the ball of radius $r$ centered at $x$ by $x+r\B=\set{y\in\R^n:\,\abs{y-x}\le r}$, where $\abs{\;\cdot\;}$ is the Euclidean norm. For a closed set $A\sub\R^n$ and $x\in\R^n$, we denote the distance from $x$ to $A$ by $\abs{x}_{A}=\inf_{y\in A}\abs{x-y}$ and $r$-neighborhood of $A$ by $A+r\B=\cup_{x\in A}(x+r\B)=\set{x\in\R^n:\,\abs{x}_A\le r}$. For a set $A\subseteq\R^n$, $\overline{A}$ denotes its closure, $\inte(A)$ denotes its interior, and $\partial A$ denotes its boundary.  
For two sets $A,B\sub\R^n$, the set difference is defined by $A\setminus B=\set{x:\,x\in A,\,x\not\in B}$.

We denote by $F:\R^m\sra\R^n $ a set-valued map.  We say a function $\alpha:\,\R_{\ge 0}\rightarrow\R_{\ge 0}$ belongs to class $\mathcal{K}$ if it is continuous, zero at zero, and strictly increasing. It is said to belong to $\mathcal{K}_{\infty}$ if it belongs to class $\mathcal{K}$ and is unbounded. A function $\beta:\,\R_{\ge 0}\times\R_{\ge 0}\rightarrow\R_{\ge 0}$ is said to belong to class $\mathcal{KL}$ if, for each $t\ge 0$, $\beta(\cdot,t)$ belongs to class $\mathcal{K}$ and, for each $s\ge 0$, $\beta(s,\cdot)$ is decreasing and satisfies $\lim_{t\rightarrow\infty}\beta(s,t)=0$.

\section{Preliminaries}\label{sec: pre}
\subsection{\bf{Hybrid systems}}
Consider a hybrid  system 
$\hbd=(C,F,D,G)$ with dynamics
\begin{subequations}\label{E: sys}
\begin{align}
& \dot{x}\in F(x),\;\;\;\;\;x\in C,\\
& x^+\in G(x), \;\;\;x\in D,
\end{align}
\end{subequations}
where $C,D\sub \R^n$ represent the flow set and the jump set respectively, and the set-valued maps $F:\R^n\sra\R^n $, $G:\R^n\sra\R^n $ represent the flow and jump maps respectively. 

Given a scalar $\delta\geq 0$, the (additive) $\delta$-perturbation of $\hbd$, denoted by $\hbd_\delta$, is described as
\begin{subequations}\label{E: sys_p}
\begin{align}
& \dot{x}\in F(x)+\delta\B,\;\;\;\;\;x\in C_\dt,\\
& x^+\in G(x)+\delta\B, \;\;\;x\in D_\dt,
\end{align}
\end{subequations}
where $C_\dt=C+\dt\B$ and $D_\dt=D+\dt\B$. 
A hybrid time domain is a subset $E\sub \R_{\geq 0}\times\N $ of the form $\cup_{j=0}^{J}[t_j,t_{j+1}]\times\set{j}$, where $J\in \N\cup\set{\infty}$ and $0=t_0\leq t_1\leq t_2\leq \dots$. Given $(t,j),(t',j')\in E$, the natural  ordering is such that $t+j\leq t'+j'$ if $t\leq t'$ and $j\leq j'$. A hybrid arc is a function $\phi:E\ra\R^n$ from a hybrid domain $E$ to $\R^n$ and, for each fixed $j$, $t\mapsto\phi(t,j)$ is locally absolutely continuous on the interval $I^j=\{t:(t,j)\in \dom(\phi)\}$.

\begin{deff}[Solution concept]
A hybrid arc $\phi$ is a solution to \eqref{E: sys} if 
\begin{enumerate}
    \item[(i)]$\phi(0,0)\in \overline{C}\cup D$; 
    \item [(ii)] for all $j\in\N$ such that $\inte(I^j)\neq \emptyset$, we have 
    \begin{equation*}
        \begin{split}
            & \phi(t,j)\in C\;\;\text{for all} \;t\in\inte(I^j),\\
            & \dot{\phi}(t,j)\in F(\phi(t,j))\;\;\text{for almost all} \;t\in I^j;
        \end{split}
    \end{equation*}
    \item [(iii)] for all $(t,j)\in\dom(\phi)$ such that $(t,j+1)\in\dom (\phi)$, we have $\phi(t,j)\in D$ and $\phi(t,j+1)\in G(\phi(t,j))$.
\end{enumerate}
Solutions to $\hbd_\dt$ for $\dt>0$ are defined in a similar way. Furthermore, a solution $\phi$ to $\hbd_\dt$ is maximal if there is no other $\psi$ to $\hbd_\dt$ such that $\phi(t,j)=\psi(t,j)$ for all $(t,j)\in\dom(\phi)$ and $\dom(\phi)$ is a proper subset of $\dom(\psi)$ \cite[Definition 2.7]{goebel2012hybrid}. 
\end{deff}

We also define the range of a solution arc $\phi$ as $\rge(\phi)=\{\phi(t,j): (t,j)\in\dom(\phi)\}$ for convenience. For $\delta\geq 0$, we denote by  $\soln_{\hbd}^{\delta}(x) $  the set of all maximal solutions 
starting from $x\in\R^n$
for a hybrid system   $\hbd_\delta$; we denote by  $\soln_{\hbd}^{\delta}(K) $  the set of all maximal solutions starting from the set $K\sub\R^n$
for $\hbd_\delta$.  

We  introduce correspondingly some notations for reachable sets of $\hbd_\delta$ from some $K\sub\R^n$. For $T\geq 0$, we define 
\begin{equation*}
    \rs^\delta_{\leq T}(K)= \set{\phi(t,j): \phi\in\soln_\hbd^\delta(K), \phi(0,0)\in K, t+j\leq T}.
\end{equation*}
The reachable sets $\rs^\delta_{< T}(K)$ and $\rs^\delta_{\geq T}(K)$ are defined in a similar way. The infinite-horizon reachable set from $K$ for $\hbd_\delta$ is $\rs^\delta(K)=\set{\phi(t,j): \phi\in\soln_\hbd^\delta(K), \;(t,j)\in\dom(\phi)}$.

\begin{ass}\label{ass: basic} We make the standing assumption that $\hbd$ should satisfy the basic conditions:
\begin{enumerate}
    \item[(C1)]  $C,D\in\R^n$ are closed;
    \item[(C2)] $F$ is outer semicontinuous, locally bounded, and convex for all $x\in C$;
    \item[(C3)] $G$ is outer semicontinuous and locally bounded  for all $x\in D$.
\end{enumerate}
\end{ass}

We also introduce the following concept of closeness of hybrid arcs that will be frequently used in proofs. 
\begin{deff}\cite[Definition 5.23]{goebel2012hybrid}
Given $\tau,\eps>0$, two hybrid arcs $\phi_1$ and $\phi_2$ are $(\tau,\eps)$-close if 
\begin{enumerate}
    \item[(a)] for all $(t,j)\in\dom(\phi_1)$ with $t+j\leq \tau$ there exists $s$ such that $(s,j)\in\dom(\phi_2)$, $\abs{t-s}<\eps$ and $\abs{\phi_1(t,j)-\phi_2(s,j)}<\eps$;
    \item[(b)] for all $(t,j)\in\dom(\phi_2)$ with $t+j\leq \tau$ there exists $s$ such that $(s,j)\in\dom(\phi_1)$, $\abs{t-s}<\eps$ and $\abs{\phi_2(t,j)-\phi_1(s,j)}<\eps$.
\end{enumerate}
\end{deff}

\subsection{\bf{Properties pertaining to stability and safety}}
In this subsection, we define properties that are related to stability and safety of solutions to \eqref{E: sys_p}. Particularly, we formally introduce the concepts of stability with safety guarantees and reach-avoid-stay type specifications.
\begin{deff}[Forward (pre-) invariance]\label{def: fw_inv}
A set $\ik\sub\R^n $ is said to be forward pre-invariant for $\hbd_\delta$ if for every $\phi\in\soln_\hbd^\delta(\ik)$, $\rge(\phi)\sub \ik$. The set $\ik$ is said to be forward invariant for $\hbd_\delta$ if for every forward complete $\phi\in\soln_\hbd^\delta(\ik)$, $\rge(\phi)\sub \ik$, i.e., $\rs^\delta(\ik)\sub \ik$.
\end{deff}
\begin{rem}
The term `pre' is in the sense that non-complete maximal solutions are not excluded. This concept allows us to describe the completeness and  dynamical behaviors separately for general maximal solutions. For future references, we only define the `pre' properties of solutions to keep the relevant definitions succinct. 
\end{rem}

We next consider stability for solutions of $\hbd_\delta$ w.r.t. a closed set.
\begin{deff}[Uniform pre-asymptotic stability]\label{def: UpAS}
A closed set $A\sub\R^n$ is said to be uniformly pre-asymptotically stable (UpAS) for $\hbd_\delta$ if the following two conditions are met:
\begin{enumerate}
    \item[(i)](uniform stability) for every $\eps>0$, there exists a $\dt_\eps>0$ such that every solution to $\hbd_\dt$  with $\abs{\phi(0,0)}_A<\dt_\eps$ satisfies $\abs{\phi(t,j)}_A<\eps$ for all $(t,j)\in\dom(\phi)$;
     \item[(ii)] (uniform pre-attractivity) there exists some $\rho>0$ such that, for every $\eps>0$, there exist some $T>0$ such that, for every solution $\phi$ to $\hbd_\dt$,
     $|\phi(t,j)|_A< \eps$ 
     whenever $|\phi(0,0)|_A<\rho$, $(t,j)\in\dom(\phi)$, and $t+j\geq T$.
\end{enumerate}
\end{deff}
\begin{deff}[Basin of pre-attraction]
If a closed set $A\sub \R^n$ is UpAS for $\hbd_\dt$, we further define the basin of pre-attraction of $A$, denoted by $\bk_\dt(A)$, as the set of initial states $x\in\R^n$ such that every solution $\phi$ to $\hbd_\dt$ with $\phi(0,0)=x$ is bounded and, if it is complete, then also converges to the set $A$, i.e., $\lim_{t+j\ra\infty}|\phi(t,j)|_A=0$.
\end{deff}

We define two closely-related properties of our interests in the following definitions. 
\begin{deff}[Stability with safety guarantee]\label{def: stable_safe}
Let $\xs_0,U,A\sub\R^n$ and $A$ is a closed set. 
We say that $\hbd_\dt$ satisfies a stability with safety guarantee specification $(\xs_0,U,A)$ if 
\begin{enumerate}
    \item[(i)](pre-asymptotic stability w.r.t. $A$) the set $A$ is UpAS for $\hbd_\dt$ and $\xs_0\sub\bk_\dt(A)$; 
    \item[(ii)](safety w.r.t. $U$) $\rge(\phi)\cap U=\emptyset$ for all $\phi\in\soln_\hbd^\dt(\xs_0)$.
\end{enumerate}
\end{deff}

\begin{deff}[Reach-avoid-stay specification]\label{def: ras}
Let $\xs_0,U,\ik\sub\R^n$. We say that $\hbd_\dt$ satisfies a reach-avoid-stay specification $(\xs_0,U,\ik)$ if 
\begin{enumerate}
    \item[(i)](reach and stay w.r.t. $\ik$) there exists some $T\geq 0$ such that $\rs^\dt_{\geq T}(\xs_0)\sub\ik$;
    \item[(ii)](safety w.r.t. $U$) $\rge(\phi)\cap U=\emptyset$ for all $\phi\in\soln_\hbd^\dt(\xs_0)$.
\end{enumerate}
\end{deff}

\section{Connection between Stability  with  Safety  Guarantees and Reach-avoid-stay  Specifications}\label{sec: connect}
\subsection{ Stability with safety implies 
reach-avoid-stay}
Throughout this section, we suppose that $\hbd_\dt$ satisfies a stability with safety
guarantee specification $(\xs_0, U, A)$ for some fixed $\dt\geq 0$.

We first show some basic properties of solutions to $\hbd_\dt$ with a fixed $\dt\geq 0$. The following proposition combines the results from  \cite[Proposition 7.4, Lemma 7.8]{goebel2012hybrid}.
\begin{prop}\label{prop: basin}
Given a hybrid system $\hbd_\dt$, let $A\sub\R^n$ be a compact set that is uniform stable for $\hbd_\dt$ (condition (i) of Definition \ref{def: UpAS} holds). 
Suppose $K\sub\bk_\dt(A)$ is compact. Then
\begin{enumerate}
    \item[(i)] $\bk_\dt(A)$ is open and $A\sub\bk_\dt(A)$;
    \item[(ii)] For every $\eps>0$, there exists some $T=T(\eps)>0$ such that 
    $
        |\phi(t,j)|_A\leq \eps
    $
       for all $\phi\in\soln_\hbd^\dt(K)$ and $(t,j)\in\dom(\phi)$ with $t+j\geq T$;
    \item[(iii)]$A\cup\rs(K) $ is a compact subset of $\bk_\dt(A)$.
\end{enumerate}
\end{prop}

By (ii) of Proposition \ref{prop: basin}, the uniform attraction for solutions to $\hbd_\dt$, the reach-avoid-stay property can be shown straightforwardly. The statement is given in the following proposition.
\begin{prop}
Let $A,\xs_0\sub\R^n$ be compact sets. Suppose that $\hbd_\dt$ satisfies a stability with safety guarantee specification $(\xs_0, U, A)$. Then for every $\eps>0$, the system $\hbd_\dt$ also satisfies the reach-avoid-stay specification $(\xs_0, U, A+\eps\B)$.
\end{prop}

\subsection{The converse side}
Throughout this subsection, we suppose that $\hbd_\dt$ satisfies a reach-avoid-stay specification $(\xs_0, U,\ik)$. We make an additional assumption on the local Lipschitz property of $F$ and $G$ to prove the converse connection. 

\begin{deff}[Locally Lipschitz set-valued maps]
Let $\mathcal{O}\sub \R^n$ be an open set. 
A set-valued flow map $F$ is also locally Lipschitz on $\mathcal{O}$; that is, for each  $x\in\mathcal{O}$, there exists a neighborhood $\mathcal{N}\subset\mathcal{O}$ of $x$ and an $L>0$ such that for any $x_1,x_2\in\mathcal{N}$, 
\begin{equation*}
    F(x_1)\sub F(x_2)+L|x_2-x_1|\B.
\end{equation*}
\end{deff}

The following fact is a direct consequence of $F$ being locally Lipschitz \cite[Lemma 9]{ teel2000smooth}.
\begin{lem}\label{lem: lip}
If $F:\mathcal{O}\sra \R^n$ is locally Lipschitz, then for any compact subset $K\sub\mathcal{O}$ there exists an $L_K>0$ such that for any $x_1,x_2\in K$, we have
\begin{equation*}
   F(x_1)\sub F(x_2)+L_K|x_2-x_1|\B.  
\end{equation*}
\end{lem}
\begin{ass}\label{ass: lip}
We assume in this section that $F$ and $G$ in \eqref{E: sys_p} are locally Lipschitz on some open subset of $C$ and $D$ respectively. 
\end{ass}

\begin{lem}[Perturbed solutions]\label{lem: construct}
For any $\dt> 0$, fix a $\dt'\in[0,\dt)$ and a $\tau>0$. Let $T\geq \tau$ and $K\sub\R^n $ be a compact set such that $K\cap(C_\dt\cup D_\dt)\neq \emptyset$. Then there exists an $r=r(K,\tau,\dt',\dt)>0$  such that for every solution $\phi\in\soln_\hbd^{\dt'}(K)$ with 
$\phi(t',j')\in K$ for all $(t',j')\in\dom(\phi)$ and $t'+j'\leq T$, and for all $x\in \phi(0,0)+r\B$, 
there exists a $\psi\in\soln_\hbd^\dt(x)$ with $\dom(\psi)=\dom(\phi)$ 
and
\begin{equation}\label{E: psi}
    \psi(t,j)=\phi(t,j),\;\;t+j=T.
\end{equation}
\end{lem}
\begin{pf}
We suppose the total number of jumps before $T$ is $N$, thereby $N\in[0,T]$.
Now, set 
\begin{equation*}
    \psi(t',j')=\phi(t',j')+\left(1-\frac{t'+j'}{T}\right)\cdot(\psi(0,0)-\phi(0,0))
\end{equation*}
for all $(t',j')\in\dom(\phi)$  and $t'+j'\in[0,T]$. Then the constructed $\psi$ satisfies \eqref{E: psi} and
\begin{equation*}
    \abs{\psi(t',j')-\phi(t',j')}\leq |\psi(0,0)-\phi(0,0)|\leq r.
\end{equation*}
 In particular, for $N>0$,  
 for each jump  $j'\in\{0,1,\cdots,N-1\}$ 
 at the end  point $t_j'$ of $I^{j'}$, we have
 \begin{equation*}
     |\psi(t_j',j'+1)-\phi(t_j',j'+1)|\leq r. 
 \end{equation*}
Similarly, suppose $N<T$ (i.e., flows exist), 
 then for any fixed $j'\in\{0,1,\cdots,N\}$ and $t'\in I^{j'} $, we have
 \begin{equation*}
     \abs{\dot{\psi}(t',j')-\dot{\phi}(t',j')}\leq \frac{1}{T}\abs{\psi(0,0)-\phi(0,0)}= \frac{r}{T}.
 \end{equation*}
 We show that this $\psi$ is a solution to $\hbd_\dt$.  On the flow set $C_\dt$, let $\wh{C}:=C_\dt\cap (K+r\B)$. Then for $x\in \wh{C}$, for all $t'\in I^{j'}$, we have 
 \begin{equation*}
 \begin{split}
     \dot{\psi}(t',j')\sub&\dot{\phi}(t',j')+\frac{r}{T}\B\sub F(\phi(t',j'))+\left(\dt'+\frac{r}{T}\right)\B\\
          \sub& F(\psi(t',j'))+\left(\dt'+\frac{r}{T}\right)\B+L_{\wh{C}}|\psi(t',j')-\phi(t',j')|\B\\
          \sub & F(\psi(t',j'))+\left(\dt'+\frac{r}{T}+rL_{\wh{C}} \right)\B\\
          \sub & F(\psi(t',j'))+\left(\dt'+\frac{r}{\tau}+rL_{\wh{C}} \right)\B,
 \end{split}
 \end{equation*}
 where the constant $L_{\wh{C}}$ in line 4 of the above inclusions is obtained by Lemma \ref{lem: lip}. 
 
 Similarly, for each end point $t'$ of some $I^{j'}$, the jump should satisfy
 \begin{equation*}
   \psi(t',j'+1)\sub G(\psi(t',j'))+ \left(\dt'+r+rL_{\wh{D}} \right)\B,  
 \end{equation*}
 where $\wh{D}:=D_\dt\cap(K+r\B)$ and $L_{\wh{D}}$ is obtained based on Lemma \ref{lem: lip} for $G$. 
 
 Now we are able to bound the above by $\dt$ such that $\psi$ is a solution  to $\hbd_\dt$. The $r$ can be selected accordingly such that all of the followings are satisfied: $r\leq \dt-\dt'$,  $\frac{r}{\tau}+\dt'+rL_{\wh{C}}\leq \dt$ and $r+\dt'+rL_{\wh{D}}\leq \dt$.
 
 Note that the above choice of $r$ should work for every possible case: the extreme cases when only flows ($N=0$) or jumps ($N=T$) happen, as well as the mixed flow/jump case ($N\in(0,T)$). 
\end{pf}
\begin{rem}
The above result shows that for any compact solution $\phi$ to $\hbd_{\dt'}$ (that exists till $t+j=\tau$), there exists a solution $\psi$ to a slightly more perturbed system $\hbd_\dt$ such that the endpoints of $\phi$ and $\psi$ are related within some time period. The construction and estimation rely on the local Lipschitz continuity of $F$ and $G$. The proof states that no matter how conservative the estimation is, we are able to find the small neighborhood with radius $r$ such that for any initial condition within $\phi(0,0)+r\B$, there exists a $\psi$ as a solution to $\hbd_\dt$ converging to $\phi$ within finite time. The requirement that $\dt'$ should be strictly less than $\dt$ is necessary. Unlike the robustness  concept given in \cite[Lemma 7.37]{goebel2012hybrid}, the perturbed systems $\hbd_\dt$  need  inflation of $C$ and $D$ whose intensities increase as $\dt$ increases. This subtle difference is in consideration of when $\phi(0,0)\in \partial(C_{\dt'}\cup D_{\dt'})$ whilst the constructed $\psi$ is still well posed.  
\end{rem}

Applying the preceding results, we show in the next proposition a hybrid-system version of \cite[Proposition 17]{meng2022smooth}. The proof is completed in Appendix. 
\begin{prop}\label{prop: connect}
Any nonempty, compact, forward pre-invariant set $A$ for $\hbd_\dt$ is UpAS for $\hbd_{\dt'}$ with any $\dt'\in[0,\dt)$.
\end{prop}

%

The existence of  such a nonempty, compact, and forward pre-invariant set for $\hbd_\dt$ is guaranteed by the following proposition. 
\begin{prop}[Existence of  compact invariant set]\label{lem: invariant}
Suppose that $\ik$ is compact and $\xs_0$ is nonempty. Suppose that $\hbd_\dt$ satisfies a reach-avoid-stay specification $(\xs_0, U,\ik)$. Then the set 
\begin{equation}
    A=\set{x\in\ik: \forall \phi\in\soln_\hbd^\dt(x), \rge(\phi)\sub\ik}
\end{equation}
is nonempty, compact, and forward pre-invariant for $\hbd_\dt$.
\end{prop}

To prove Proposition \ref{lem: invariant}, we need the following lemma.

\begin{lem}\label{lem: initial}
Let the hypothesis in Proposition \ref{lem: invariant} be satisfied. Let $x_0\in A\sub\ik$ be fixed. Then for every $\eps>0$ and $\tau\geq 0$, there exists $\kappa>0$ such that, for every solution $\phi\in\soln_\hbd^\dt(x_0+\kappa\B)$ there exists a solution $\psi\in\soln_\hbd^\dt(x_0)$ such that $\phi$ and $\psi$ are $(\tau,\eps)$-close. 
\end{lem}

\begin{pf}
Without loss of generality, we  assume that $\eps$ and $\tau$ are arbitrarily small such that $x_0+\kappa\B\sub\ik$. By the reach-avoid-stay property of $\hbd_\dt$, the solution starting from $\ik$ should be either eventually bounded or complete. Suppose the statement were to fail, then for some arbitrarily small $\eps$ and $\tau$, for each $n\in\N$ and $\phi_n\in\soln_\hbd^\dt(x_0+1/n\B)$ with $x_0+1/n\B\sub\ik$, there exists no solution $\psi\in\soln_\hbd^\dt(x_0)$ is $(\tau,\eps)$-close to $\phi_n$. However, $\{\phi_n\}_n$ is  (locally eventually) bounded. By \cite[Theorem 6.1]{goebel2012hybrid}, we can extract a subsequence, still denoted by $\{\phi_n\}_n$, that is graphically convergent to some $\phi\in\soln_\hbd^\dt(x)$. This implies that there exists some sufficiently large $n$ such that $\phi_n$ and $\phi$ are $(\tau,\eps)$-close as a consequence of graphical convergence of bounded sequences, which leads to a contradiction. 
\end{pf}
\textbf{Proof of Proposition \ref{lem: invariant}}: 
By Definition \ref{def: ras}, there exists some $T\geq 0$ such that 
$\rs_{\geq T}^\dt(\xs_0)\sub\ik$. It is easy to verify that for any $x\in \rs_{\geq T}^\dt(\xs_0)\sub\ik$, $x\in\ik$, and for all $\phi\in\soln_\hbd^\dt(x)$, we have $\phi(t,j)\in\rs_{\geq T}^\dt(\xs_0) $ for all $(t,j)\in\dom(\phi)$. This shows that $A$ is nonempty. The forward pre-invariance is verified by setting $T=0$ and $\xs_0\sub\ik$.

It suffices to show the closedness of $A$, which will imply the compactness due to the boundedness assumption on $\ik$. Pick any sequence $\set{x_n}_n\sub A$ such that $x_n\ra x$. 
Then for all $\phi_n\in\soln_\hbd^\dt(x_n)$, we have $\rge(\phi_n)\sub\ik$ and hence bounded. Suppose that $x\notin A$, then there exists a  $\phi\in\soln_\hbd^\dt(x)$ such that $\phi(t,j)\notin\ik$ for some $(t,j)\in\dom(\phi)$. Select sufficiently small $\tau$ and $\eps$ as in Lemma \ref{lem: initial}, then there exists some $\kappa=\kappa(\tau,\eps)>0$ and sufficiently large $n\in\N$ (with $x\in x_n+\kappa\B$)  such that $\phi\in\soln_\hbd^\dt(x_n+\kappa\B)$. However, by Lemma \ref{lem: initial}, there exists a solution $\phi_n\in\soln_\hbd^\dt(x_n)$ such that $\phi$ and $\phi_n$ are arbitrarily $(\tau,\eps)$-close. Since $\R^n\setminus{\ik}$ is open, we have $\phi_n(s,j)\in\R^n\setminus{\ik}$ for some $s+j\leq \tau$ with $|t-s|<\eps$. This contradicts the forward pre-invariance property of $\phi_n$. Hence, $x\in A$ and $A$ depicts compactness. \hfill\BlackBox\\[2mm]

Combining Proposition \ref{prop: connect} and Proposition \ref{lem: invariant}, we immediately obtain the following connection between reach-avoid-stay specifications and stability with safety guarantees.
\begin{prop}\label{prop: connect_2}
Let Assumption \ref{ass: lip} be satisfied. Let $\xs_0$ be a nonempty set and $\ik$ be a compact set. Suppose that $\hbd_\dt$ satisfies a reach-avoid-stay specification $(\xs_0, U,\ik)$. Then there exists a compact set $A\sub\ik$ such that any less perturbed system $\hbd_{\dt'}$ with $\dt'\in[0,\dt)$ satisfies the stability with safety guarantee specification $(\xs_0, U,A)$.
\end{prop}

\section{Lyapunov-Barrier Functions for
Reach-Avoid-Stay Specifications}\label{sec: converse}
In this section, we provide Lyapunov-barrier characterizations for stability with safety guarantee specifications as well as reach-avoid-stay specification. In particular, such characterizations are (robustly) complete in the sense that  Lyapunov-barrier functions exist based on the specified dynamical behavior of the solutions. 

We focus on the  Lyapunov-barrier functions for  stability with safety guarantees, the result for reach-avoid-stay comes after the connection given in Proposition \ref{prop: connect_2}. Since the uniform asymptotic stability or reachability of a closed set strongly relies the set initial conditions and may not hold globally, the point-to-set distance becomes inadequate on the topology of basin of pre-attraction. We use a proper indicator to characterize the distance instead. 
\begin{deff}[Proper indicator]
Let $A\sub\R^n$ be a compact set and $\ok\sub\R^n$ be an open set containing $A$. A continuous function $\omega: \ok\ra\R_{\geq 0}$ is said to be a proper indicator for $A$ on $\ok$ if
\begin{enumerate}
    \item[(i)] $\omega(x)=0\iff x\in A$;
    \item[(ii)] $\lim_{n\ra \infty}\omega(x_n)=\infty$ for any sequence $\set{x_n}_n\sub \ok$ such that either $x_n\ra p\in\partial \ok$ or $|x_n|\ra \infty$ when $n\ra \infty$.
\end{enumerate}
\end{deff}

We also refine the basin of pre-attraction and introduce the following constraint set.
\begin{deff}[Basin  of  pre-attraction with safety]\label{def: BA_safe}
The extracted basin  of  pre-attraction with safety is a subset of $\bk_\dt(A)$ defined by
\begin{equation}
    \wh{\bk}_\dt(A):=\{x\in\bk_\dt(A): \forall \phi\in\soln_\hbd^\dt(x),\;\;\rge(\phi)\cap U = \emptyset\},
\end{equation}
where $A,U$ are given in the specification $(\xs_0,U,A)$.
\end{deff}

\ymmark{The following lemma verifies basic topological properties of $\wh{\bk}_\dt(A)$, which will be later utilized for the Lyapunov-barrier function theorems. }
\begin{lem}\label{lem: basin_open}
Suppose that $A$ is compact, $U$ is closed, and $A\cap U=\emptyset$. If $\hbd_\dt$ satisfies the stability with safety guarantee specification $(\xs_0,U,A)$, then
\begin{enumerate}
    \item[(i)] $\xs_0\subseteq\wh{\bk}_\dt(A)\sub \bk_\dt(A)$;
    \item[(ii)]$\wh{\bk}_\dt(A)$ is open and forward pre-invariant.
\end{enumerate}
\end{lem}
\begin{pf}
The first claim can be easily verified by Definition \ref{def: stable_safe} and \ref{def: BA_safe}. The forward pre-invariance of (ii) comes after Definition \ref{def: fw_inv}. It suffices to show that $ \wh{\bk}_\dt(A)$ is also open given $U$ is closed. 

Suppose the opposite, then pick any $x\in \wh{\bk}_\dt(A)$. There exists a sequence of points $\set{x_n}_n\sub\bk_\dt\setminus\wh{\bk}_\dt(A)$ with $x_n\ra x$. Note that for each $n$, there exists $\phi_n\in\soln_\hbd^\dt(x_n)$ that is either bounded or complete with convergence to $A$. Either way, we have $\phi_n$ bounded for each $n$ by Proposition \ref{prop: basin}. By \cite[Theorem 6.1]{goebel2012hybrid}, we can extract a subsequence, still denoted by $\set{\phi_n}_n$, that is graphically convergent. By Assumption \ref{ass: basic}, the limit satisfies $\phi\in\soln_\hbd^\dt(x)$, and is again either bounded or complete with convergence  to $A$.

Suppose the graphical limit $\phi$ is complete. Let $r>0$ be given by condition (ii) from Definition \ref{def: UpAS}. Pick $\eps<r$ such that $(A+\eps\B)\cap U=\emptyset$. Choose $\dt_\eps\leq \eps$ according to condition (i) of Definition \ref{def: UpAS} by the uniform stability. Then there exists a $T=T(\dt_\eps)>0$ such that 
$|\wh{\phi}(t,j)|_A\leq \dt_\eps$
for all $\wh{\phi}\in\soln_\hbd^\dt(x)$ and $t+j\geq \dom(\wh{\phi})$ with $t+j\geq T$. Note that the reachable set $\rs_{\leq T}^\dt(x)$ is compact and we can select some $\eps'<\eps/2$ such that $(\rs_{\leq T}^\dt(x)+\eps'\B)\cap U=\emptyset$. Let $\eps''=\min\{\eps', \dt_\eps\}$. As the consequence of graphical convergence, there exists a sufficiently large $n\in\N$ such that for each $\tau>0$, the solution $\phi_n$ and $\phi$ are $(\tau,\eps'')$-close.
No matter $\tau\geq T$ or $\tau<T$, we can verify by the choice of $\eps''$ and $n$ that $\phi_n(t,j)\in (\rs^\dt(x)+\eps'\B)\cup (A+\eps\B)$ for all $(t,j)\in\dom(\phi_n)$, which contradicts the property of $\phi_n\in\soln_\hbd^\dt(x_n)$ with $x_n\in\bk_\dt(A)\setminus\wh{\bk}_\dt(A)$. For the case that $\phi$ is bounded, we can proceed and show the contradiction by a similar way. Combining the above, we have $\wh{\bk}_\dt(A)$ is open.
\end{pf}
 
The following result shows that,  based on the reachable region of solutions with stability and safety guarantees, a single Lyapunov-like function exists which is also effective as a barrier function to guarantee the safety. 

\begin{thm}\label{thm: lya_stability}
Suppose that $A$ is compact, $U$ is closed, and
$A\cap U=\emptyset$. Then the following two statements are equivalent:
\begin{enumerate}
    \item[(i)] $\hbd_\dt$ satisfies the stability with safety guarantee specification $(\xs_0,U,A)$.
    \item[(ii)] There exists an open and forward pre-invariant set $\ok$ such that $(A\cup \xs_0)\sub \ok$ and $\ok\cap U=\emptyset$, a smooth function $V: \ok\ra\R_{\geq 0}$ and class $\mathcal{K}_\infty$ functions $\alpha_1,\alpha_2$ such that, 
    \begin{equation}
        \alpha_1(\omega(x))\leq V(x)\leq \alpha_2(\omega(x)),\;\;\forall x\in (C_\dt\cup D_\dt\cup G(D_\dt))\cap \ok,
    \end{equation}
    \begin{equation}
        \nabla V(x) \cdot f\leq -V(x), \;\;\forall x\in C_\dt\cap \ok, \;f\in F(x)+\dt\B,
    \end{equation}
    and 
    \begin{equation}
        V(g)\leq V(x)/e,\;\;\forall x\in D_\dt\cap \ok, \;g\in G(x)+\dt\B,
    \end{equation}
    where $\omega$ is a proper indicator for $A$ on $\ok$.
\end{enumerate}
\end{thm}
\begin{pf}
(i) $\implies$ (ii): By \cite[Corollary 7.33]{goebel2012hybrid}, the existence of smooth Lyapunov-like function holds on forward pre-invariant open subsets of $\bk_\dt(A)$. The result follows immediately by Lemma \ref{lem: basin_open} for $\ok=\wh{\bk}_\dt(A)$.

    (ii) $\implies$ (i): By a standard Lyapunov argument for hybrid systems, we are able to show that for any $\phi$ to $\hbd_\dt$ with $\phi(0,0)\in \ok$, 
    \begin{equation*}
        V(\phi(t,j))\leq V(\phi(0,0))e^{-(t+j)/3},
    \end{equation*}
    which implies the forward pre-invariance of $\ok$ and the UpAS property of any $\phi$ starting within $\ok$. Now that $\ok\cap U=\emptyset$, for any $\phi\in\soln_\hbd^\dt(\xs_0)$ with $\xs_0\sub\ok$, we have $\rge(\phi)\cap U=\emptyset$, which implies the safety. 
\end{pf}

\begin{rem}
Note that to use \cite[Corollary 7.33]{goebel2012hybrid}, we need to verify that the UpAS holds for the following system of the following modifications:
\begin{subequations}\label{E: sys_r}
\begin{align}
& \dot{x}\in F_\rho(x),\;\;\;\;\;x\in (C_\dt)_\rho,\\
& x^+\in G_\rho(x), \;\;\;x\in (D_\dt)_\rho,
\end{align}
\end{subequations}
where 
\begin{equation*}
    (C_\dt)_\rho=\set{x\in\R^n: (x+\rho(x)\B\cap C_\dt\neq\emptyset},
\end{equation*}
\begin{equation*}
    F_\rho(x)=\overline{\operatorname{con}}F((x+\rho(x)\B)\cap C_\dt)+(\rho(x)+\dt)\B,
\end{equation*}
\begin{equation*}
    (D_\dt)_\rho=\set{x\in\R^n: (x+\rho(x)\B\cap D_\dt\neq\emptyset},
\end{equation*}
and $G_\rho(x)=\set{v\in \R^n: v\in g+(\rho(g)+\dt)\B}$, for $g\in G(x+\rho(x)\B)\cap D_\dt$.  However, this is guaranteed by Assumption \ref{ass: basic} for each $\dt\geq 0$. 

The proof for \cite[Corollary 7.33]{goebel2012hybrid} relies on the conversion from UpAS to the $\mathcal{KL}$ pre-asymptotic stability\footnote{The definition is omitted. We kindly refer readers to \cite[Definition 31]{liu2020smooth} for more details.}  based on changing the basin of attraction. In contrast, for robust systems driven by differential equations with $F$ being a single-valued function, \cite[Proposition 10]{liu2020smooth} provides a more explicit way of constructing the $\mathcal{KL}$-function $\beta$ such that $\omega(\phi(t))\leq\beta(\omega(x),t)$ for all solutions with initial conditions $x\in\ok$. The argument is based on the reachable set properties as well as a connection between the solution-to-set distance and the proper indicator. The extension of the explicit construction from \cite[Proposition 10]{liu2020smooth} to the general hybrid systems should be similar. 
\end{rem}

By the connection obtained in Section \ref{sec: connect}, the following results follow.

\begin{cor}
Suppose that $A,\xs_0$ are compact sets and (ii) of Theorem \ref{thm: lya_stability} holds. Then $\hbd_\dt$ satisfies  the  reach-avoid-stay  specification $(\xs_0,U,A+\eps\B)$ for any $\eps>0$.
\end{cor}

\begin{thm}
Let Assumption \ref{ass: lip} be satisfied. 
Suppose that $\ik$ is compact, $U$ is closed and $\ik\cap U=\emptyset$. Suppose that $\hbd_{\dt'}$ satisfies the reach-avoid-stay specification $(\xs_0, U,\ik)$.  Then there exists a compact set $A\sub\ik$ such that (ii) of Theorem \ref{thm: lya_stability} holds for any $\hbd_\dt$ with $\dt\in[0,\dt')$. 
\end{thm}

The above constructed Lyapunov functions for the stability with safety guarantees and reach-avoid-stay specifications  play an implicit role as a barrier function. As though it appears more succinct to have both Lyapunov and barrier functions merged into one single Lyapunov function, it may not be easy to design in practice.  We provide an equivalent characterization as in Theorem \ref{thm: lya_stability} with separate Lyapunov-barrier functions to complete this section. The separate Lyapunov-barrier functions for reach-avoid-stay specifications are omitted due to the similarity.

\begin{thm}\label{thm: lya-barrier}
Suppose that $A$ is compact, $U$ is closed and $A\cap U=\emptyset$. If there exists an open set $\ok$ such that $A\cup \xs_0\sub \ok$, a smooth function $V: \ok\ra\R_{\geq 0}$ satisfying 
\begin{enumerate}
    \item[(i)] there exist  $\alpha_1,\alpha_2\in\mathcal{K}$ and a continuous positive function $\varrho$ such that
    \begin{equation}
    \begin{aligned}
        &\alpha_1(|x|_A)\leq V(x)\leq \alpha_2(|x|_A),\;\;\forall x\in (C_\dt\cup D_\dt\cup G(D_\dt))\cap\ok,
        \end{aligned}
    \end{equation}
    \begin{equation}
         \nabla V(x)\cdot f\leq -\varrho(|x|_A),\;\;\forall x\in C_\dt,\;f\in F(x)+\dt\B,
    \end{equation}
    \begin{equation}
         V(g)-V(x)\leq -\varrho(|x|_A),\;\;\forall x\in D_\dt,\;g\in G(x)+\dt\B;
    \end{equation}

\end{enumerate}
and $B: \R^n\ra \R$ that is smooth in $C\cap\ok$ such that
\begin{enumerate}
        \item[(ii)] the set $S:=\set{x\in\R^n: B(x)\geq 0}\sub\ok$ and $\xs_0\sub S$;
    \item[(iii)] $x\in U$ implies $B(x)<0$;
    \item[(iv)]
    \begin{equation}
        \nabla B(x)\cdot f\geq  0,\;\;\forall x\in C_\dt\cap \ok, \;f\in F(x)+\dt\B,
    \end{equation}
    \begin{equation}
       B(g)-B(x)\geq 0,\;\;\forall x\in D_\dt\cap  \ok, \;g\in G(x)+\dt\B,
    \end{equation}
\end{enumerate}
then $\hbd_\dt$ satisfies the stability with safety guarantee specification $(\xs_0,U,A)$.

If $\xs_0$ is also compact, then the converse side also holds, i.e., the existence of smooth $V$ and $B$ with conditions (i)-(iv) is necessary for $\hbd_\dt$ satisfying the stability with safety guarantee specification $(\xs_0,U,A)$.
\end{thm}
\begin{pf}
For the sufficient part, we start with the stability. Note that by a standard Lyapunov argument for hybrid systems, the existence of $V$ with condition (i) is sufficient to guarantee the UpAS property of $A$ for $\hbd_\dt$. The condition (ii)-(iv)  of $B$ intends to separate the set $S$ and $U$, such that all the solutions starting from $S$ will stay within it. 

We show formally the mechanism of $B$. Suppose the opposite, then there exists some $x\in S$, a solution $\phi\in\soln_\hbd^\dt(x)$, and some hybrid time such that $B(\phi(t,j))<0$. Then we are able to define a finite time
\begin{equation*}
    \tau:=\sup\{t+j\geq 0: \phi(t,j)\in S\}.
\end{equation*}
As a consequence, there exist $t'$ and $j'$ with $t'+j'=\tau$ such that $B(\phi(t',j'))=0$. Note that since the $\phi$ exists till $B(\phi(t,j))<0$. For a small perturbation $\eps$ of $t'$, $\phi$ is either a flow such that $t'+\eps\in I^{j'}$, or triggers a jump such that $t'+\eps\in I^{j'+1}$. For the first case, since $\ok$ is open, for arbitrary $\eps>0$, we still have $\phi(t'',j')\in C_\dt\cap\ok$ for almost $t''\in[t',t'+\eps]$. Thus, we have $\dot{B}(\phi(t'',j'))=\nabla B(\phi(t'',j))\cdot f\geq 0$ for all $f\in F(x)+\dt\B$ and almost all $t''\in[t',t'+\eps]$, which means $t''+j'\geq \tau$ but $\phi(t'',j')\in S$. This contradicts the definition of $\tau$. For the second case, by a similar argument,  we have $B(\phi(t',j'+1))\geq 0$, which also leads to a contradiction.  Hence $S$ is forward pre-invariant. This verifies condition (ii) of Definition \ref{def: stable_safe}.

A direct consequence of $B$ is that, for any $\xs_0\sub\ok$, any solution  $\phi\in\soln_\hbd^\dt(\xs_0)$ stays within $S$ and hence $\ok$. Condition (i) guarantees that $\phi\in\bk_\dt(A)$. This verifies condition (i) of Definition \ref{def: stable_safe}. 

Now we show the converse side. The existence of $V$ is derived from the existing converse Lyapunov theorems given the UpAS property of the compact set $A$. The construction of $B$ is also based on the compactness of $A$ and $\xs_0$. 
Let $c=\sup_{x\in \xs_0} V(x)$. Then, $\xs_0\sub K:=\{x\in\ok: V(x)\leq c\}$. \ymmark{We aim to unite locally smooth functions of the form $B_i(x)=c-V(x)$ into a globally smooth function $B$. The rest of the proof falls in  the exact procedure as in the proof of \cite[Proposition 12]{meng2022smooth}. We choose $c_2>c_1>c$ and let $\ok_i=\{x\in\ok: V(x)<c_i\}$ for $i\in\{1,2\}$. Then $\{\ok_1, \ok_2\setminus K, \R^n\setminus\overline{\ok_1}\}$ is an open cover of $\R^n$. On these converings, let $B_1(x)=c-V(x)$ for $x\in\ok_1$, $B_2(x)=c-V(x)$ for $x\in \ok_2\setminus K$, and $B_3(x)=c-c_1$ for $x\in \R^n\setminus\overline{\ok_1}$. The above functions are defined to be $0$ elsewhere. Let ${\psi_1, \psi_2, \psi_3}$ be a smooth partition of unity \cite[p.43]{lee2013smooth} subordinate to $\{\ok_1, \ok_2\setminus K, \R^n\setminus\overline{\ok_1}\}$. Then, $B(x)=\sum_{i=1}^3\psi_i(x)B_i(x)$ for $x\in\R^n$ can be verified to satisfy condition (ii)-(iv).}
\end{pf}

\begin{rem}
For the sufficient part, condition (iv) intends to regulate the safe direction of  solutions on the entire open set $\ok$ rather than directly on its subset $S:=\set{x\in\R^n: B(x)\geq 0}$. This condition seems stronger but also necessary. Imposing condition (iv) on $S$ can only guarantee the invariance of the interior of $S$, for counters examples see e.g. \cite[Remark 4]{wang2021safety}.
\end{rem}

\begin{rem}
The necessity of the existence of  $B$ given the satisfaction of a stability with safety guarantee specification $(\xs_0,U,A)$ is majorly based on the compactness of reachable set from $\xs_0$, which does not intersect with $U$ by the \textit{a priori} safety assumption. This compactness in turn makes it easier to construct $B$ via the Lyapunov function $V$. The latest converse barrier function theorems by \cite{ghanbarpour2021barrier,maghenem2022converse,maghenem2022strong}, however, does not depend on the boundedness of $\R^n\setminus U$, and hence does not require the compactness of $\xs_0$. The construction of barrier function given the safety of robust systems relies on the  `time-to-impact' function $B_\rs$ w.r.t. the reachable set $\rs$ of robust systems. This $B_\rs$ turns out to satisfy condition (ii) and (iii) of Theorem \ref{thm: lya-barrier} as well as a weaker certificate\footnote{Note that we have adapted the notion to be consistent with the notation used in this paper.}
\begin{equation*}
    \nabla B\cdot f>0,\;\;\forall x\in\partial S, f\in F(x)+\dt\B.
\end{equation*}
The results in  \cite{ghanbarpour2021barrier,maghenem2022converse,maghenem2022strong} are of great theoretical interest in the sense of constructing barrier functions with relaxed topological requirement and barrier conditions. Their construction relies on the seemingly less intuitive time-to-impact functions, whereas our approach aims to unify Lyapunov and barrier functions. It would be interesting to investigate in future work whether the Lyapunov-barrier approach can be extended to handle unbounded reachable sets. 
\end{rem}

\section{Examples}\label{sec: example}
We provide two examples in this subsection to validate our results.

\subsection{Bouncing ball}

As a classical mechanical system with impulse-momentum change, the model of bouncing balls is frequently used to illustrate dynamical behaviors of hybrid systems. 
We  model a tennis ball as a point-mass, and consider dropping it from a fixed height with a constant horizontal speed. The vertical direction is forced by the gravity. As the ball hits the horizontal surface, the instantaneous vertical velocity is reversed  with a small dissipation of energy.

To describe the hybrid dynamics, the state of the ball is given as 
\begin{equation*}
    \xk=(x,y,z)^T\in\R^3,
\end{equation*}
where $x$ is the horizontal position, $y$ represents the height above the surface, and $z$ is the vertical velocity. The flow set is given  as
$
    C=\{\xk: y>0\;\text{or}\; y=0, z\geq 0\}.
$
We consider the flow\footnote{As in \cite{goebel2012hybrid}, it is natural to set $f(\zero)=\zero$ regardless of the noncontinuity at $\zero$.} as
\begin{equation*}
    \dot{\xk}=f(\xk):=\begin{bmatrix}
    1\\
    z\\
    -a
    \end{bmatrix},\;\;\xk\in C,
\end{equation*}
where $a=9.8$ is the acceleration due to the gravity. The jump set captures the domain when the vertical velocity flips the sign, which is geiven as 
$
    D=\{\xk: y=0\;\text{and}\; z<0\}.
$
The jump dynamic is such that
\begin{equation*}
    \xk^+=g(\xk):=\begin{bmatrix}
    x\\
    0\\
    -\varsigma z
    \end{bmatrix},\;\;\xk\in D,
\end{equation*}
for some dissipation coefficient $\varsigma\in(0,1)$. 

Now fix $\varsigma=0.8$, let the initial condition to be $\xk(0)=(0,9,0.8)^T\in C$,  and model a block (the unsafe set) by 
$
    U:=\{\xk: y>10\}.
$
We consider the target set as 
$
   \ik:=\{\xk: y\in[0, 0.1]\}. 
$

\begin{figure}
    \centering
    \includegraphics[scale=0.50]{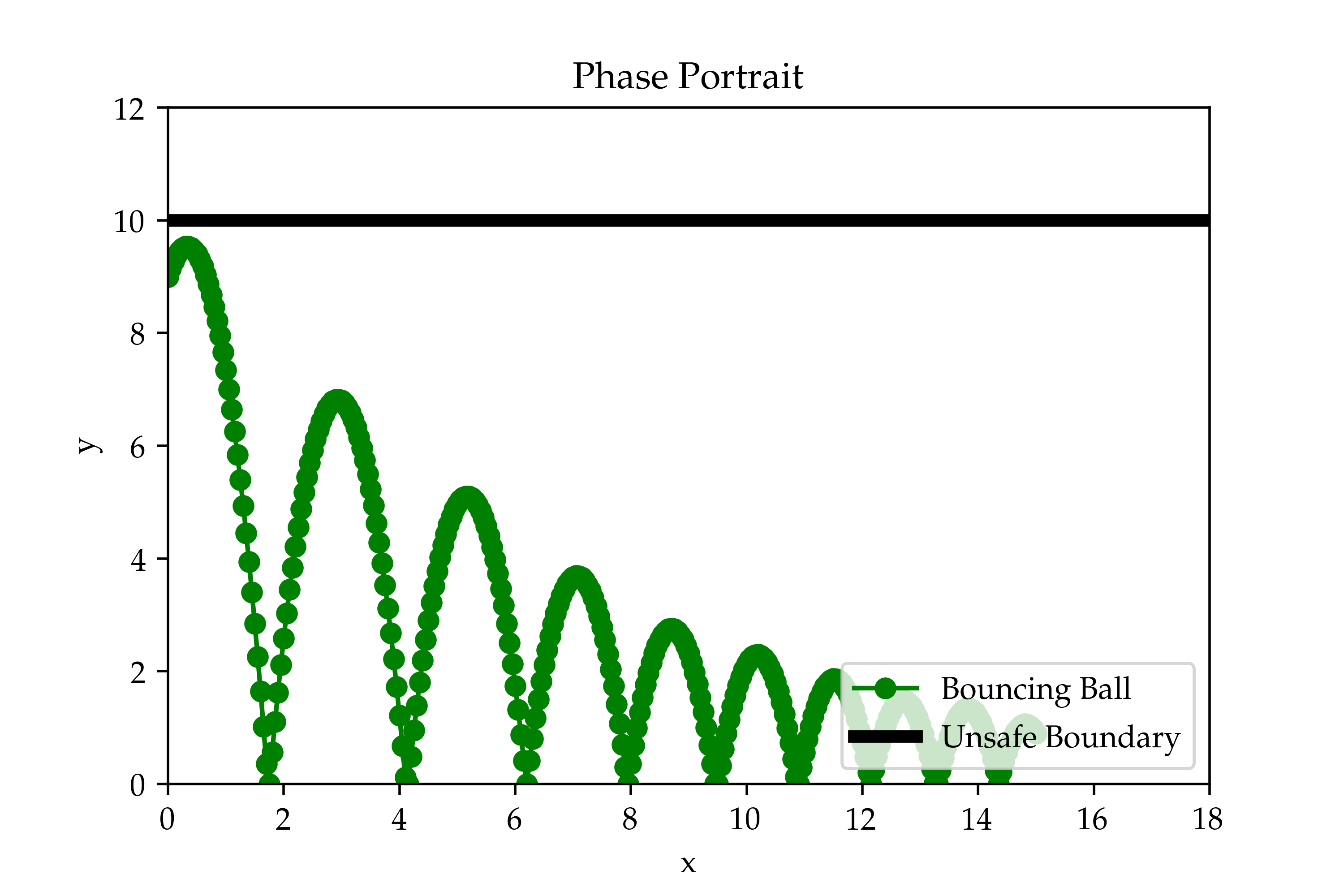}
    \caption{Snapshot of bouncing ball in $xy$-plane.}
    \label{fig: bouncing_2d}
\end{figure}

\begin{figure}
    \centering
    \includegraphics[scale=0.6]{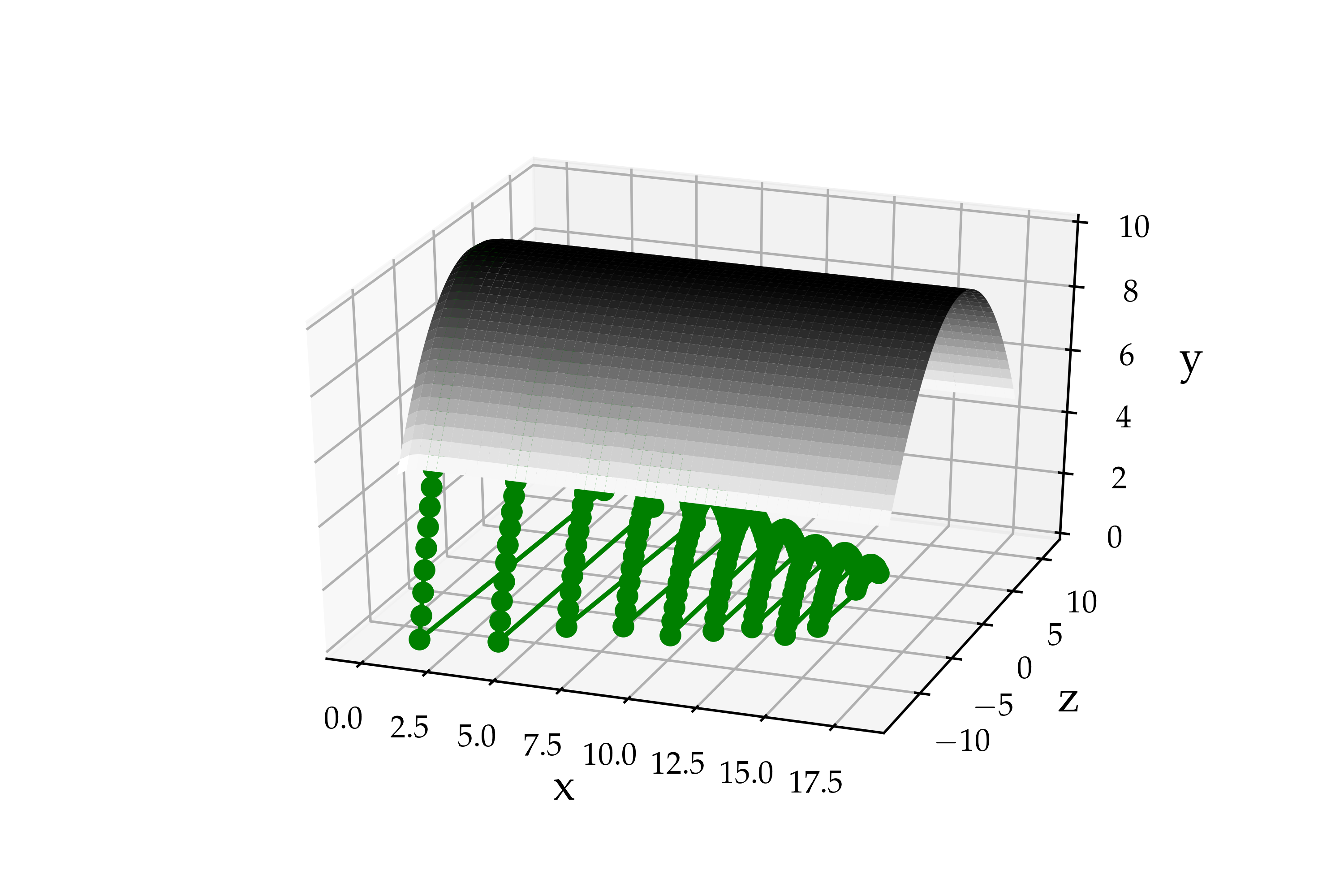}
     \includegraphics[scale=0.6]{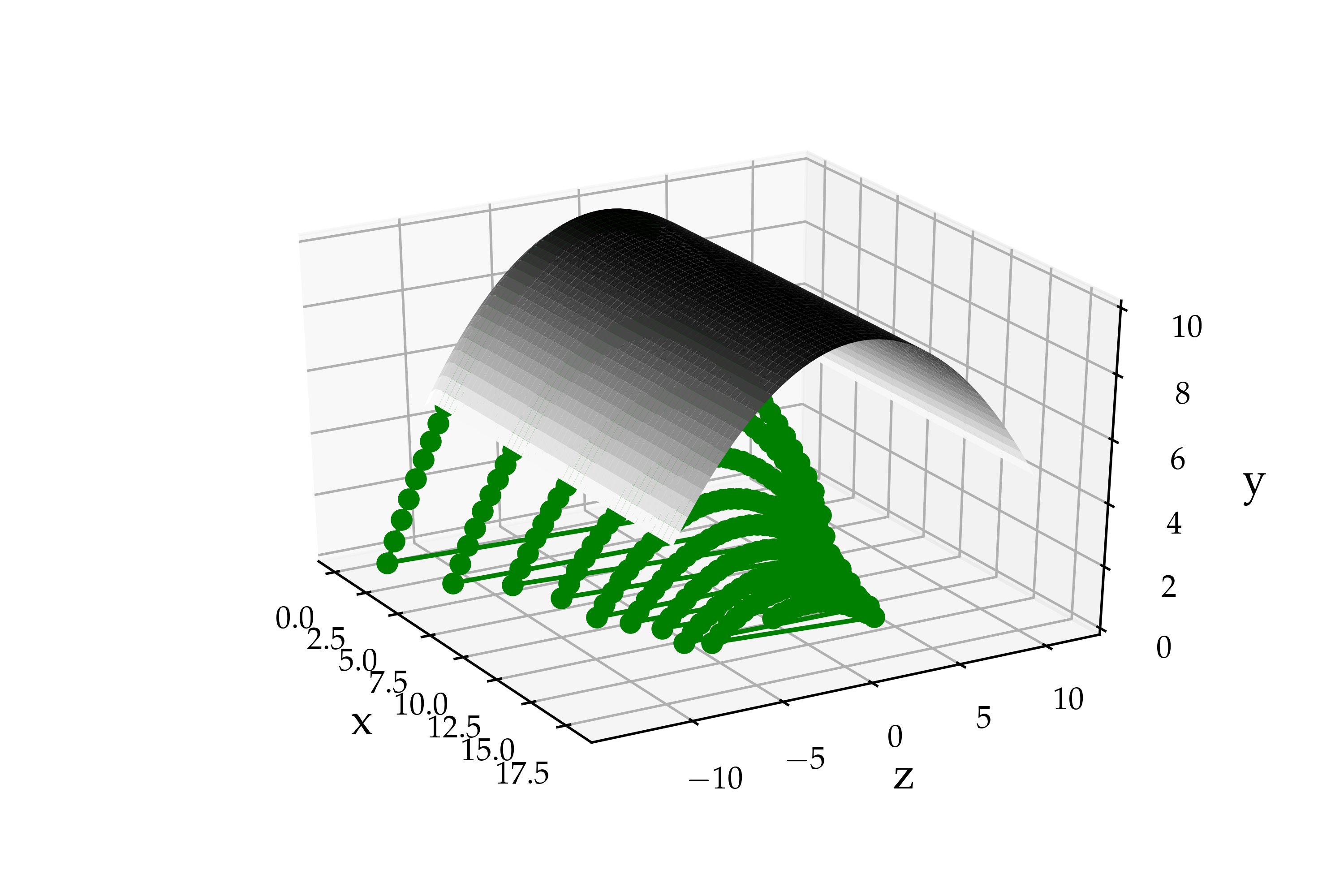}
    \caption{Solution of bouncing ball (the green dots) and barrier function $B$ (the grey surface) in $xyz$-plane.}
    \label{fig: bouncing_3d}
\end{figure}

It can be shown analytically and  numerically (see Figure \ref{fig: bouncing_2d}) that the system satisfies the reach-avoid-stay specification $(\{\xk(0)\}, U,\ik)$. 
We now  validate the existence of  Lyapunov-barrier functions $(V,B)$. 
Note that, the function 
\begin{equation*}
    V(\xk)=\left(1+\frac{1-\varsigma^2}{\pi(1+\varsigma^2)}\arctan(z)\right)\cdot\left(\frac{1}{2}z^2+ay\right)
\end{equation*}
has been verified 
in \cite[Example 3.19]{goebel2012hybrid} to be a valid Lyapunov function w.r.t. the set $\{\xk: y=0\}$. It suffices to find valid candidate of the barrier function $B$. Consider a sigmoid function
\begin{equation*}
   \sigma(x)=\frac{1}{1+\exp(-5x)} 
\end{equation*}
and let
\begin{equation*}
    B(\xk)=\frac{1}{2}\sigma(x)-y-\frac{1}{2a}z^2+9.5.
\end{equation*}
It can be verified that $B(\xk)<0$ for $\xk\in U$, 
\begin{equation*}
    \nabla B(\xk)\cdot f(\xk)=\frac{1}{2}\sigma(x)(1-\sigma(x))-z+z\geq 0,\;\;\xk\in C,
\end{equation*}
and 
\begin{equation*}
    B(g(\xk))-B(\xk)=\frac{1}{2a}(z^2-\varsigma^2z^2)>0,\;\;\xk\in D,
\end{equation*}
Therefore, the above $B(\xk)$ is a valid barrier function. 
Let $S:=\{\xk: B(\xk)\geq 0\}$ be the set as in Theorem \ref{thm: lya-barrier}, then it is also clear that $\xk(0)\in S$. The evolution of the state $\xk\in\R^3$ as well as the graph of barrier function $B$ are provided in Figure \ref{fig: bouncing_3d}. It can be seen that for all $t\geq 0$, $\xk(t)\in S$ and therefore remain safety w.r.t. $U$.

\subsection{Sample-and-hold control}
In this example, we revisit the case study in \cite[Section V]{meng2021control} from the perspective of hybrid system. 

We consider the reduced Moore-Greitzer ODE model with an additive control input $[v,0]^T$ and no extra perturbations:
\begin{equation}\label{E:MG}
    \dfrac{d}{dt}\begin{bmatrix}
\Phi(t)\\ \Psi(t)
\end{bmatrix}=\begin{bmatrix}
\frac{1}{\mathfrak{l}_c}(\psi_c-\Psi(t))\\
\frac{1}{16\mathfrak{l}_c}\left(\Phi(t)-\gamma\sqrt{\Psi(t)}\right)
\end{bmatrix}+\begin{bmatrix}
v(t)\\0
\end{bmatrix},
\end{equation}
where $\psi_c=\mathfrak{a}+\iota\cdot[1+\frac{3}{2}\left(\frac{\Phi}{\Theta}-1\right)-\frac{1}{2}\left(\frac{\Phi}{\Theta}-1\right)^3]$, and $\gamma$ is the throttle coefficient. The other parameters are as follows:
\begin{equation*}
    \mathfrak{l}_c=8,\;\iota=0.18,\;\Theta=0.25,\;\mathfrak{a}=1.67\iota.
\end{equation*}
As
 $\gamma$ decreases, surge instability occurs
and generates a pumping oscillation (Hopf-bifurcation) that
can cause flameout and engine damage. The physical meanings of variables, parameters and the
description of this model can be found in \cite[Section V]{meng2021control}.

\begin{prob}\label{P: case_study_det}
We aim to manipulate the throttle coefficient $\gamma$ and $v$ simultaneously such that the state $(\Phi,\Psi)$ are regulated to satisfy reach-avoid-stay specification $(\xs_0,U,\Gamma)$. We require that $\gamma:\R_{\geq 0}\rightarrow [0.5,1]$ is time-varied with $\gamma(0)\in[0.62,0.66]$ and  $|\gamma(t+\tau)-\gamma(t)|\leq 0.01\tau$ for any $\tau>0$.  We define $\xs_0=\{(\Phi_e(\gamma(0)),\Psi_e(\gamma(0)))\}$
(i.e. a sub-region of stable equilibrium points, also see \cite[Fig. 1]{meng2021control}); $\Gamma$ to be the ball that centered at $\zeta=(0.4519,0.6513)$ with radius $r=0.003$, i.e. $\Gamma=\zeta+r\mathcal{B}$; $U=\{(x,y):x\in(0.497,0.503),y\in(0.650,0.656)\}$. We set 
$v(t)\in\mathcal{U}=[-0.05,0.05]\cap \R$ for all $t$. \pfbox
\end{prob}
Note that, the system \eqref{E:MG} can be rewritten as 
\begin{equation}\label{E:sys5}
    \dfrac{d}{dt}\begin{bmatrix}
\Phi(t)\\ \Psi(t)
\end{bmatrix}=\begin{bmatrix}
\frac{1}{\mathfrak{l}_c}(\psi_c-\Psi(t))\\
\frac{1}{16\mathfrak{l}_c}\Phi(t)
\end{bmatrix}+\begin{bmatrix}
v(t)\\-\gamma(t)\sqrt{\Psi(t)}
\end{bmatrix}.
\end{equation}
Let $\uf(t)=[v(t),\gamma(t)]^T$ and 
\begin{equation*}
    g(x(t))=\begin{bmatrix}
    1 & 0\\
    0 & -\sqrt{\Psi(t)}
    \end{bmatrix},
\end{equation*}
then \eqref{E:sys5} fits in the general form of flow dynamics. We make a little abuse of notation here to define $\mathcal{U}=[-0.05,0.05]\times [0.5,1]$.

Due to the issues of less frequent and inaccurate state measurement,  the errors inject into the 
closed-loop control system  and may cause unsatisfactory performance. We are motivated by the above issues to convert the system \eqref{E:sys5} into a hybrid system. We then apply the conditions of the Lyapunov-barrier functions in Theorem \ref{thm: lya-barrier} to synthesize valid sample-and-hold controllers aiming at fulfilling the task in Problem \ref{P: case_study_det}.

For simplicity, we define the state $\xk=(\Phi,\Psi)^T$ and the control input $\uf=(v,\gamma)$. We write
\begin{equation*}
    \tilde{f}(\xk,u)= \begin{bmatrix}
\frac{1}{\mathfrak{l}_c}(\psi_c-\Psi(t))\\
\frac{1}{16\mathfrak{l}_c}\Phi(t)
\end{bmatrix}+\begin{bmatrix}
v(t)\\-\gamma(t)\sqrt{\Psi(t)}
\end{bmatrix},
\end{equation*}
which is the r.h.s. of \eqref{E:sys5}. 

We then introduce $\tau\in [0,0.5)$ as a timer variable and set
\begin{equation*}
    \zk=\begin{bmatrix}
\xk\\
\uf\\
\tau\end{bmatrix}\in\R^2\times\R^2\times\R.\end{equation*}
The flow set of the sample-and-hold hybrid system is given as 
$C=\R^2\times\R^2\times[0,0.5)$,  
on which the flow dynamics are such that 
\begin{equation*}\dot{\zk}=f(\zk)=\begin{bmatrix}
\tilde{f}(\xk,\uf)\\
0\\
1
\end{bmatrix},\;\;\zk\in C.
\end{equation*}

$$$$
The jumps happen when the timer is up to $0.5$ and decisions are made, i.e., 
$D=\R^2\times\R^2\times\{0.5\}$
and 
\begin{equation}
    \zk^+=g(\zk)=\begin{bmatrix}
\xk\\
\kappa(\xk)\\
0
\end{bmatrix},\;\;\zk\in D.
\end{equation}

We choose Lyapunov-barrier functions $(V,B)$ to be such that 
\begin{itemize}
    \item[(1)] $V(x)=\|x-\zeta\|^2$  for all $x\in D\setminus \{\zeta\}$,
    \item[(2)] and $B(x)=-\log\left(\frac{h(x)}{1+h(x)}\right)$ for $h(x)=\|x-(0.500,0.653)\|_{\infty}-r$.
\end{itemize}
We then define a set of high-gain robust control policy as follows: 
\begin{equation}\label{E: jump_MG}
\begin{split}
     \mathfrak{K}(x)=&\{u\in\uu: L_fV(x)+L_gV(x)u+V(x)\leq -\varsigma\;\text{and}\; L_fB(x)+L_gB(x)u\geq \varsigma,\;\forall x\in D\},
\end{split}
\end{equation}
where $\varsigma>0$ is intended to compensate the error that is generated when the same control input is imposed on the flow dynamics. Since $L_fV$, $L_gV$, $L_fB$ and $L_gB$ are locally Lipschitz continuous, and the quantity $|\xk(\tau)-\xk(0)|$ for any fixed $\xk(0)\in D$ is locally bounded given $\xk(\tau)=\xk(0)+\int_0^\tau\tilde{f}(\xk(s),u)ds$ and $\tau\in[0,0.05)$,  it is clear that both 
\begin{equation*}
\begin{split}
    &|L_fV(\xk(\tau))+L_gV(\xk(\tau))u+V(\xk(\tau))-L_fV(\xk(0))+L_gV(\xk(0))u+V(\xk(0))|
\end{split}
\end{equation*}
and 
\begin{equation*}
    |L_fB(\xk(\tau))+L_gB(\xk(\tau))u-L_fB(\xk(0))+L_gB(\xk(0))u|
\end{equation*}
are locally bounded. We  set the bound for both of the above quantities to be $\varsigma$, such that for any $\uf=\kappa(\xk)\in\mathfrak{K}(\xk)$ and $\xk\in D$, then
the Lyapunov-barrier conditions are satisfied in the flow, i.e., 
\begin{equation*}
    L_fV(\xk)+L_gV(\xk)u+V(\xk)\leq 0,\;\forall \xk\in C,\end{equation*} and \begin{equation*} L_fB(\xk)+L_gB(\xk)\uf\geq 0,\;\forall \xk\in C.
\end{equation*}
This in turn guarantees the reach-avoid-stay property of the system. 

\begin{figure}
    \centering
    \includegraphics[scale=0.50]{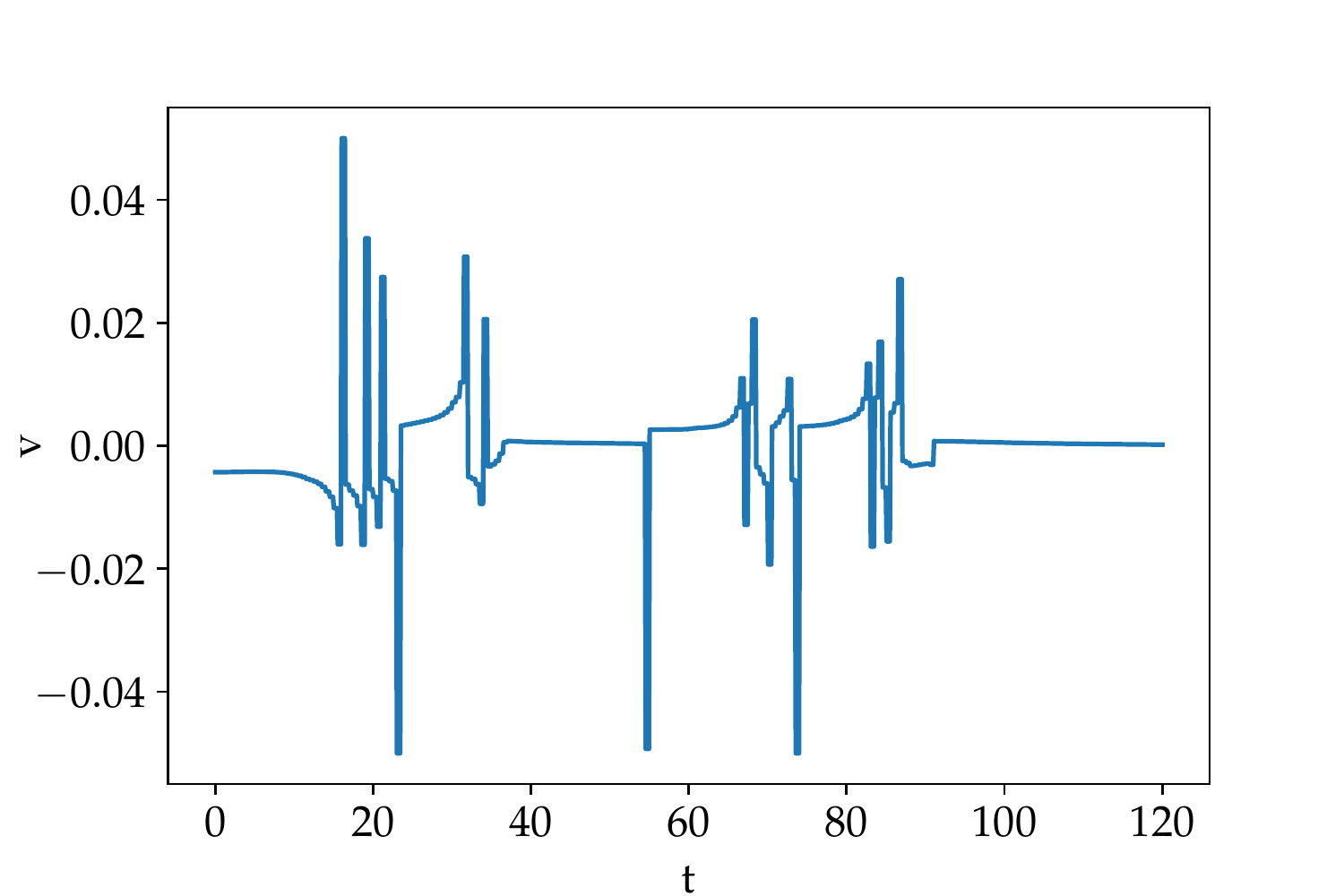}
     \includegraphics[scale=0.50]{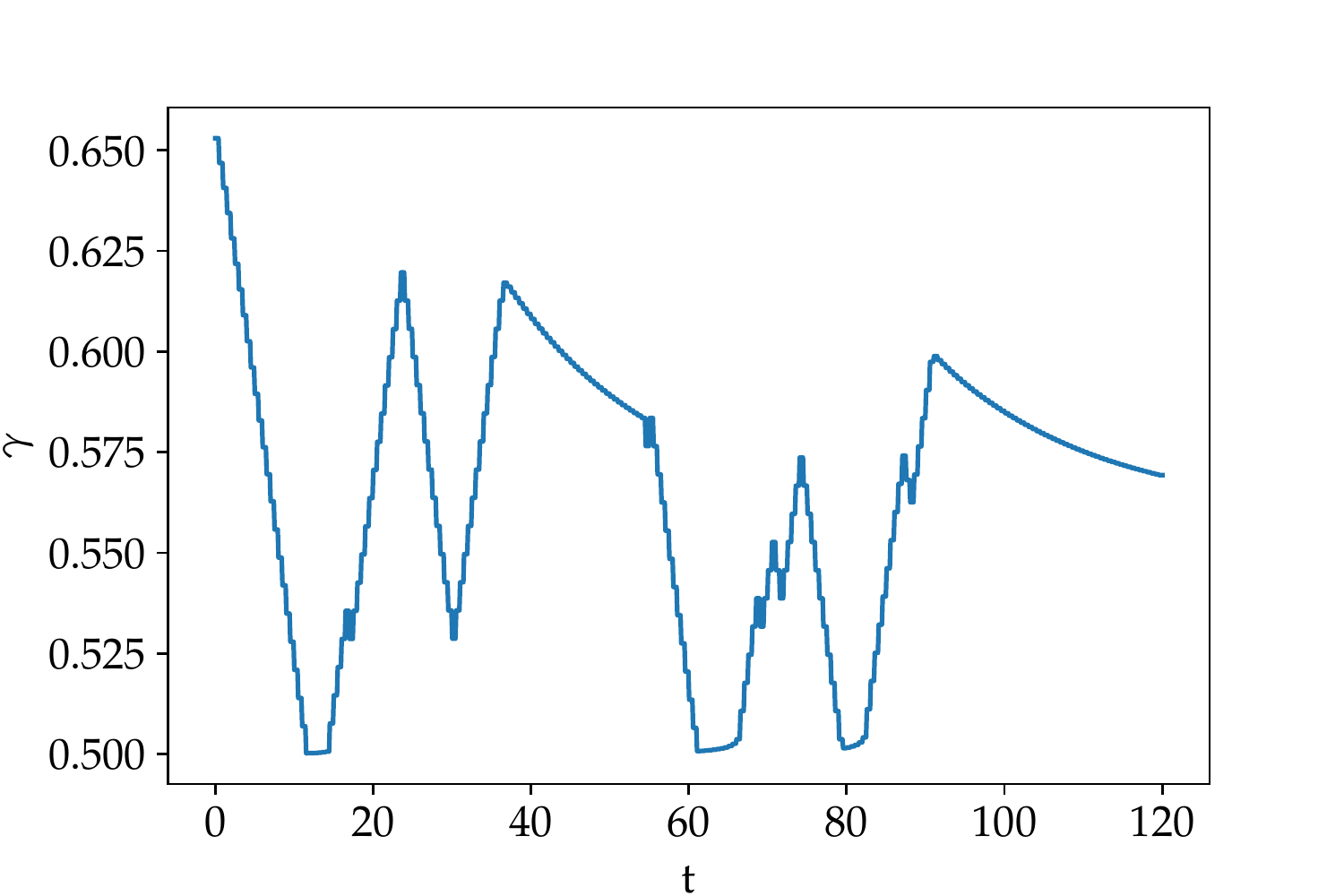}
    \caption{Control signal $v$ and $\gamma$ solved by the quadratic programming for the hybrid conversion of \eqref{E:sys5}.}
    \label{fig: hbd_MG}
\end{figure}

\begin{figure}
    \centering
    \includegraphics[scale=0.5]{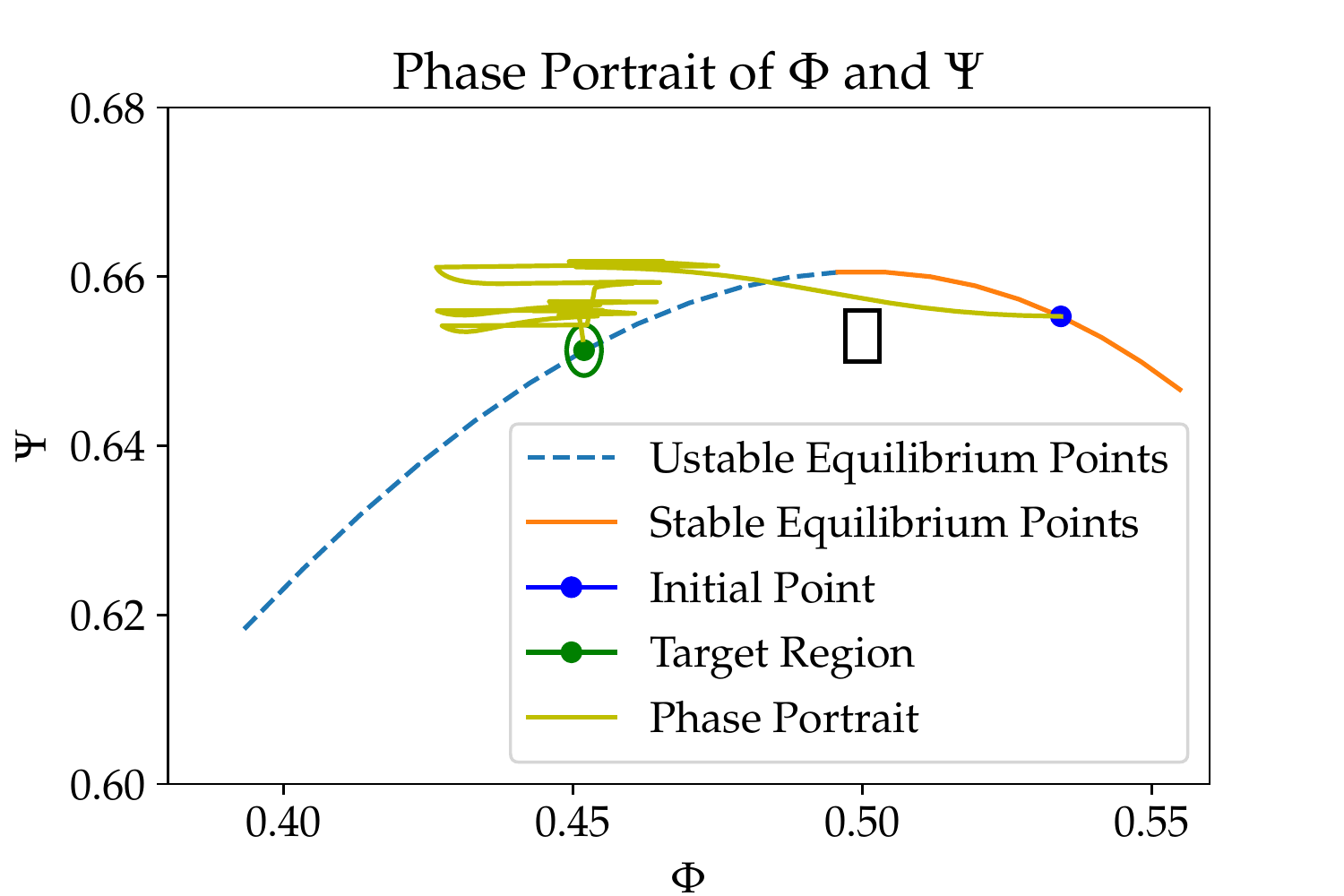}
    \includegraphics[scale=0.5]{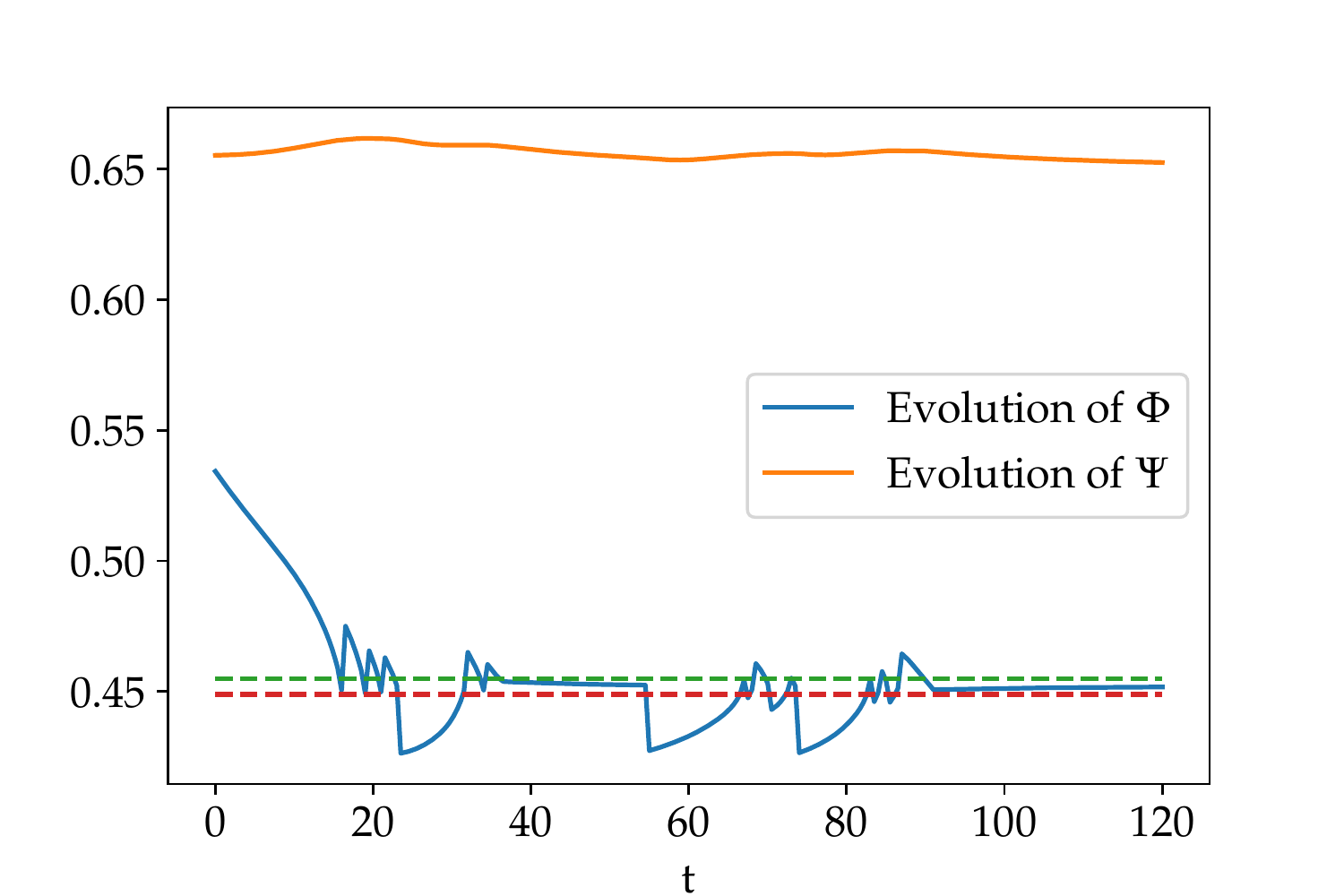}
    \caption{Phase portraits and evolution of $\Phi$ and $\Psi$ given the sample-and-hold control signal $v$ and $\gamma$.}
    \label{fig: hbd_MG_soln}
\end{figure}

In the numerical experiments, we empirically set $\varsigma=0.07$ and impose an extra set of constraints for the control input as 
\begin{equation}\label{E: constraint_MG}
\begin{split}
    \mathcal{M}& :=\{\gamma\in[0.5,1]: \gamma(0)\in[0.62,0.66],\;|\gamma(t+\tau)-\gamma(t)|\leq 0.01\tau,\;\forall \tau>0\}.
\end{split}
\end{equation}
The $\uf$ is selected such that the control effort
\begin{equation}\label{E: cost}
    |\uf(t)|^2+\frac{2\uf(t)}{\mathfrak{l}_c}(\psi_c-\Psi(t))+\left(\frac{1}{4\mathfrak{l}_c}(\Phi(t)-\gamma(t))\sqrt{\Psi(t)}\right)^2
\end{equation}
is minimized for each $\zk\in D$. The results are shown in Figure \ref{fig: hbd_MG} and \ref{fig: hbd_MG_soln}. It can be seen that given the robust correction in the decision making, the high-gain controller can still fulfill the goal apart from the intensively fluctuating trajectory along the $\Phi$ direction in between $t=20$ and $t=80$. 

\section{Conclusion}\label{sec: conclusion}
In this paper, we showed that under mild conditions, the connection can be made, via a robustness argument,  between stability with safety guarantees and reach-avoid-stay specifications for robust hybrid systems driven by differential and difference inclusions. We further extended the Lyapunov-barrier function theorems to robust hybrid systems that  satisfy stability with safety guarantees and reach-avoid-stay specifications. Under the concept of solutions to hybrid systems, as well as natural requirements on the compactness of target set and the set of initial conditions, we showed that the existence of Lyapunov-barrier functions is  necessary to the two specifications of our interests. 

In particular, two versions of Lyapunov-barrier functions are provided and discussed in Section \ref{sec: converse}. The single Lyapunov-like function that also plays a role as a barrier function appears to be succinct but may be  difficult to design in practice. We discussed the construction of separate Lyapunov-barrier functions in particular for stability with safety guarantee specifications. \ymmark{Numerical simulations are provided to validate our results.} It is interesting to compare with the latest converse theorems on barrier functions for systems driven by differential inclusions \cite{ghanbarpour2021barrier,maghenem2022converse,maghenem2022strong}. The mentioned reference provides a  construction  based on time-to-impact functions w.r.t.  the robust reachable sets, which in turn to be qualified as  barrier functions. In terms of safety, this method allows us to consider relaxations for the topological set-ups of target sets and the set of initial conditions, however, at the cost of not relying on the existing converse Lyapunov theorems and possibly sacrificing the smoothness of barrier functions. 

Nonetheless, we provide the theoretical results of converse theorems and hope to shed some light on designing powerful algorithms, e.g. learning techniques \cite{ravanbakhsh2017learning, zhao2020synthesizing, berkenkamp2016safe},  for constructing Lyapunov-barrier
functions with respect to stability with safety guarantee or reach-avoid-stay specifications. 

For future work, it would be an interesting application in using the current Lypunov-barrier characterizations as abstractions for linear temporal logic specifications for hybrid systems. It is also valuable to design algorithms based on this theoretical work to construct Lyapunov-like functions that improves the estimated region of attractions. 
Another interesting direction is to
explore the stochastic analogue of Lyapunov-barrier function theorems, particularly the converse theorems,  for stochastic hybrid systems.





\bibliographystyle{elsarticle-num-names}   
\bibliography{root} 

\appendix
\section{Proof of Proposition \ref{prop: connect}}
\textbf{Proof:}
The proof is similar to  \cite[Proposition 16]{liu2020smooth}. 
We need to verify condition (i) and (ii) of Definition \ref{def: UpAS}.

(i) For any $\eps>0$, there exists a $\tau>0$ such that for any  solution $\phi$ of $\hbd_{\dt'}$ with $\phi(0,0)\in A+\frac{\eps}{2}\B$ and $|\phi(t,j)|_A\geq \eps$, we have $t+j\geq \tau$. Now let $r=r(A+\eps\B,\tau,\dt',dt)$ from Lemma \ref{lem: construct} and set $\dt_\eps=\min(r,\eps/2)$. Let $\phi$ be any  solution of $\hbd_{\dt'}$ with $|\phi(0,0)|_A<\dt_\eps$. It suffices to show that $|\phi(t,j)|_A<\eps$ for all $(t,j)\in\dom(\phi)$. Suppose the opposite, then there exists $(t,j)\in\dom(\phi)$ and $t+j\geq \tau>0$ such that $|\phi(t,j)|_A\geq \eps$. Since $A$ is compact, by the selection of $\dt_\eps$, we can always find an $x\in A\cap (\phi(0,0)+r\B)$. By Lemma \ref{lem: construct}, there exists a $\psi\in\soln_\hbd^{\dt'}(x)$ with $\dom(\psi)=\dom(\phi)$ and $\psi(t,j)=\phi(t,j)\notin A$, which violates the invariance assumption of $A$ for $\hbd_\dt$. 

(ii)For any $\ep>0$, by the proof of (i), the uniform stability for $\hbd_\dt$, we can find a $\dt_\ep$ such that for any solution $\phi$ of $\hbd_\dt'$, we have $|\phi(t,j)|_A<\ep$ for all $(t,j)\in\dom(\phi)$ whenever $|\phi(0,0)|_A<\dt_\ep$. Let $r=r(A+\ep\B, \tau, \dt',\dt)$ from Lemma \ref{lem: construct} for some well-posed $\tau>0$. Choose $\rho\in(0,r)$, we now show that, for any solution $\phi$ of $\hbd_{\dt'}$, $\phi(t,j)\in A$ for all $(t,j)\in\dom(\phi)$ and $t+j\geq \tau$ whenever $|\phi(0,0)|_A<\rho$. Suppose the opposite. Then there exists some $(t,j)\in\dom(\phi)$ and $t+j\geq \tau$ such that $\phi(t,j)\notin A$. Note that we can always find an $x\in A\cap(\phi(0,0)+r\B)$ based on the choice of $\rho$. By Lemma \ref{lem: construct}, there exists a $\psi\in\soln_\hbd^{\dt'}(x)$ with $\dom(\psi)=\dom(\phi)$ and $\psi(t,j)=\phi(t,j)\notin A$, which violates the invariance assumption of $A$ for $\hbd_\dt$.  \hfill\BlackBox\\[2mm]
\end{document}